\documentclass[11pt,a4paper]{article}
\usepackage[T1]{fontenc}
\usepackage{lmodern,amsmath,amssymb,amsthm}
\usepackage[hmargin=25mm,top=40mm,bottom=25mm]{geometry}
\usepackage[hidelinks,hyperfootnotes=false]{hyperref}

\newtheorem{theorem}{Theorem}[section]
\newtheorem{corollary}[theorem]{Corollary}
\newtheorem{lemma}[theorem]{Lemma}
\newtheorem{proposition}[theorem]{Proposition}
\theoremstyle{remark}
\newtheorem{remark}[theorem]{Remark}
\DeclareMathOperator{\Stab}{Stab}
\DeclareMathOperator{\conv}{conv}
\DeclareMathOperator{\diam}{diam}
\DeclareMathOperator{\Id}{Id}
\DeclareMathOperator{\Gr}{Gr}
\DeclareMathOperator{\diag}{diag}
\DeclareMathOperator*{\argmin}{argmin}
\newcommand{\R}{\mathbb R}
\newcommand{\N}{\mathbb N}
\newcommand{\K}{\mathcal K}
\newcommand{\C}{\mathcal C}
\newcommand{\Sn}{\mathbb{S}^{n-1}}
\newcommand{\dH}{d_{\mathrm H}}
\allowdisplaybreaks[1]
\title{Iterated Steiner symmetrization of convex bodies always converges in shape}
\author{Jacopo Ulivelli}
\date{}

\begin{document}
\maketitle

\begin{abstract}
We prove that every sequence of iterated Steiner symmetrals of a convex body converges in shape. That is, every such sequence can be corrected by a sequence of rotations to achieve convergence in Hausdorff distance. Moreover, for a fixed corresponding family of subspaces, the same rotations can be used for every initial convex body, and their successive increments can be chosen to be square-summable.
\end{abstract}

\begingroup
\renewcommand{\thefootnote}{}
\endgroup
\setcounter{footnote}{0}

\section{Introduction}

The story of Steiner symmetrization started almost two centuries ago with the work of Jakob Steiner on the isoperimetric inequality for convex bodies (that is, non-empty compact and convex subsets of $\R^n$), based on the following remarkable property: Steiner symmetrization preserves the Lebesgue measure of a convex body while decreasing its surface area. Moreover, given a convex body $K\subset\R^n$, it is always possible to find a sequence of hyperplanes such that the corresponding iterated Steiner symmetrals converge to a centred Euclidean ball in the Hausdorff metric. For more details, we refer the reader to the exposition by Bianchi and Gronchi~\cite{exposition}.

The first general result in $\R^n$ is attributed to Gross~\cite{gross}, who proved that for every convex body $K$ there exists a symmetrization process converging to a centred Euclidean ball. Decades later, Mani-Levitska~\cite{mani} proved a much stronger statement: Choosing the sequence of hyperplanes independently and uniformly at random, and applying the corresponding symmetrization process to a convex body with non-empty interior, the corresponding sequence of symmetrals will converge almost surely to a centred Euclidean ball. See Burchard and Fortier \cite{burchardfortier} and Coupier and Davydov \cite{coupier} for further results on convergence of random Steiner symmetrizations. Even if these results do not leave too much space to non-converging symmetrization processes, a counterexample was eventually found in the work of Bianchi, Burchard, Gronchi, and Vol\v{c}i\v{c}~\cite{bbgv}. In the same work, the notion of convergence in shape was introduced in order to fix the counterexample: A sequence $(L_m)$ of convex bodies is said to converge in shape if there exists a sequence of rotations $(R_m)\subseteq SO(n)$ such that the adjusted sequence $(R_mL_m)$ converges in the Hausdorff metric. In~\cite[Theorem~2.2]{bbgv}, it is then proved that a class of Steiner symmetrization processes which includes the aforementioned counterexample does indeed always converge in shape. In that case, the successive angles between the symmetrizing directions are square-summable, and the same rotations work for all initial compact sets.

It is natural to ask whether every sequence of Steiner symmetrals converges in shape, a problem posed in Question~6.2 of~\cite{bbgv}. In this paper, we give a positive answer in the class of convex bodies. Additionally, inspired by a conjecture of Bianchi and Gronchi~\cite{private}, we show that every non-convergent sequence nevertheless admits correcting rotations with square-summable successive increments, a property shared by the counterexample~\cite{bbgv}. In fact, for a given sequence of hyperplanes, these rotations can be chosen once for all initial convex bodies.

This work is completed by some consequences and extensions to more general symmetrizations, a framework introduced by Bianchi, Gardner, and Gronchi~\cite{symgeom,symproc}. In particular, our results extend to Minkowski symmetrization which, in turn, allows us to analogously treat compact subsets of $\R^n$, following a procedure from~\cite{klain}.

\medskip
\paragraph{Main result.}
We start by introducing some notation. We denote by $\K^n$ the family of convex bodies in $\R^n$ and by $\Gr_{n,n-1}$ the linear hyperplanes in $\R^n$. For $H\in\Gr_{n,n-1}$ and $K\in\K^n$ the corresponding Steiner symmetral is $S_HK$ (see Section~\ref{sec:shape} for the precise definition of Steiner symmetrization). For a sequence $(H_m)\subset\Gr_{n,n-1}$, the corresponding Steiner symmetrization process produces, for every $K_0\in\K^n$, the iterated sequence of symmetrals
\[
K_m=S_{H_m}K_{m-1},\qquad m\in\N.
\]
Convergence of convex bodies is always understood with respect to the Hausdorff distance, denoted by $\dH$. Moreover, we denote by $d_{SO(n)}$ the usual distance on $SO(n)$ induced by its standard Riemannian structure (see Section~\ref{sec:square} for the details). Our main result reads as follows.

\begin{theorem}\label{thm:main}
Consider a sequence of hyperplanes $(H_m)\subset\Gr_{n,n-1}$. Then there exists a sequence of rotations $(R_m)\subset SO(n)$ such that, for every $K_0\in\K^n$, the sequence of symmetrals $K_m=S_{H_m}K_{m-1}$, $m\in\N$, satisfies
\[
\lim_{m\to\infty}\dH(R_mK_m,L)=0
\]
for some $L\in\K^n$. The rotations can be chosen with square-summable successive increments. That is, \[ \sum_{m\in\N}d_{SO(n)}(R_m,R_{m+1})^2<\infty.\] Additionally, if $(K_m)$ does not converge for some $K_0\in\K^n$, then these rotations satisfy
\begin{equation}\label{eq:dichotomy}
\sum_{m\in\N}d_{SO(n)}(R_m,R_{m+1})=\infty,
\qquad
\sum_{m\in\N}d_{SO(n)}(R_m,R_{m+1})^2<\infty.
\end{equation}
\end{theorem}

Notice that if $(K_m)$ converges for every $K_0\in\K^n$, one can trivially choose $R_m=\Id$, the identity matrix, for every $m\in\N$. Therefore, Theorem~\ref{thm:main} can be understood as the following dichotomy: Either the symmetrization process converges for every initial body, or all its sequences of symmetrals converge in shape with the same rotations $(R_m)$ satisfying~\eqref{eq:dichotomy}. While the proof of convergence in shape for a single body uses basic theory of convex bodies, the summability statement and the choice of a common sequence of rotations draw on additional elementary tools from the theory of Lie groups and spherical harmonics. For this reason, we correspondingly split the proof of Theorem~\ref{thm:main} between Section~\ref{sec:shape} and Section~\ref{sec:square}.

\medskip
\paragraph{Structure of the paper.}
In Section~\ref{sec:shape}, we introduce further notation and definitions, upon which we build the proof of convergence in shape for a single body. Moreover, we show how~\eqref{eq:dichotomy} adapts to the case of ellipses in the plane, providing a preliminary intuition of the convergence mechanism. Section~\ref{sec:square} is dedicated to proving that condition in complete generality and to choosing the rotations independently of the initial body. Finally, in Section~\ref{sec:extensions}, we show how the core constructions apply to the more general framework of symmetrizations introduced in~\cite{symgeom,symproc}.

\section{Convergence in shape}\label{sec:shape}

We recall that by $\K^n$ we denote the family of non-empty compact and convex subsets of $\R^n$, known as convex bodies. More generally, $\C^n$ denotes the non-empty compact subsets of $\R^n$. For $K\in\K^n$, its support function is $h_K(u)=\sup_{x\in K}(x\cdot u)$, $u\in\Sn$, where $\Sn$ is the Euclidean unit sphere of $\R^n$ and $x\cdot y$ is the standard scalar product in $\R^n$. Both $\K^n$ and $\C^n$ are endowed with the topology induced by the Hausdorff metric. See Schneider's monograph~\cite{schneider} for further details on the theory of convex bodies.

Consider $H\in\Gr_{n,n-1}$ and $K\in\K^n$. The Steiner symmetral $S_HK$ is obtained by replacing every non-empty interval $K\cap(x+H^\perp)$, $x\in H$, by the (possibly degenerate) interval of the same length centred at $x$ and orthogonal to $H$. The resulting set is again compact and convex. We only use Steiner symmetrization on $\K^n$ below. The final compact-set extension concerns Minkowski symmetrization.

For $C,D\in\C^n$, their Minkowski sum is the set $C+D=\{x+y:x\in C,\ y\in D\}$. Moreover, for $H\in\Gr_{n,n-1}$ we denote by $\rho_H$ the linear reflection in $\R^n$ with respect to $H$. A subset of $\R^n$ is $H$-symmetric if it is invariant for $\rho_H$. For $C\in\C^n$ and $H\in\Gr_{n,n-1}$, the Minkowski symmetral of $C$ with respect to $H$ is the compact set $M_HC=(C+\rho_H C)/2$. If, additionally, $C$ is convex, then so is $M_HC$. When restricted to $\K^n$, Steiner and Minkowski symmetrizations satisfy the inclusion
\begin{equation}\label{eq:inclusion}
S_HK\subseteq M_HK\qquad\text{for every }K\in\K^n,\ H\in\Gr_{n,n-1}.
\end{equation}
See, for example,~\cite[Lemma~2.1]{coupier}. In Section~\ref{sec:extensions}, we will explain how the generalization of this inclusion in the framework of~\cite{symgeom,symproc} allows us to extend our results.

The aim of this section is to prove convergence in shape for a fixed initial body, which we state here for the convenience of the reader.

\begin{theorem}\label{thm:one}
Consider $K\in\K^n$, a sequence $(H_m)\subset\Gr_{n,n-1}$, and the corresponding sequence of Steiner symmetrals
\[
K_0=K,\qquad K_m=S_{H_m}K_{m-1},\quad m\in\N.
\]
Then, there exists a sequence of rotations $(R_m)\subset SO(n)$ and $L\in\K^n$ (both allowed to depend on $K$ and $(H_m)$) such that
\[
\lim_{m\to\infty}\dH(R_mK_m,L)=0.
\]
\end{theorem}

The result above has the following immediate consequence.

\begin{corollary}\label{cor:one}
All subsequential limits of $K_m$ are the same up to a rotation in $SO(n)$.
\end{corollary}

\subsection{Proof of Theorem~\ref{thm:one}}

We start with some preliminary observations. First, we recall that given $K_0\in\K^n$ and $(H_m)\subset\Gr_{n,n-1}$, choosing $a>0$ large enough so that $K_0\subseteq aB^n$, where $B^n$ is the Euclidean unit ball, the sets $K_m=S_{H_m}K_{m-1}$ satisfy $K_m\subseteq aB^n$ for every $m\in\N$ since Steiner symmetrization preserves inclusions. For the remainder of this subsection the set $K_0$, the sequence $(H_m)$, and the corresponding symmetrization process will be fixed.

In terms of support functions, writing $f_m=h_{K_m}+a$ for $m\geq 0$, we infer
\[
0\leq f_m\leq 2a
\]
for every $m\geq 0$. Denoting by $O(n)$ the orthogonal linear group of $\R^n$, for $Q\in O(n)$ and $f\in C(\Sn)$ we write $\mathcal{U}_Qf(u)=f(Q^{-1}u)$. Consider now the transformations
\[
T_m=\frac12\bigl(\Id+\mathcal{U}_{\rho_{H_m}}\bigr).
\]
They are self-adjoint orthogonal projections on $L^2(\Sn)$ and hence $1$-Lipschitz. Moreover,~\eqref{eq:inclusion} implies
\begin{equation}\label{eq:order}
0\leq f_m\leq T_mf_{m-1}.
\end{equation}
Indeed, by the definition of Minkowski symmetrization,
\[
T_m f_{m-1}=h_{M_{H_m}K_{m-1}}+a.
\]
Thus~\eqref{eq:order} follows from~\eqref{eq:inclusion}. Denote by $\langle\cdot,\cdot\rangle_2$ the standard scalar product on $L^2(\Sn)$ and by $\|\cdot\|_2$ the corresponding norm. From~\eqref{eq:order} and since $T_m$ is $1$-Lipschitz, we infer 
\[
0\leq\|f_m\|_2\leq\|f_{m-1}\|_2.
\]
In particular, the sequence $(\|f_m\|_2)$ converges to a value $s\geq0$ since it is non-increasing.

Next, we introduce the following pseudo-distance, which will be the key ingredient of the proof of Theorem~\ref{thm:one}. For $K,K'\in\K^n$, consider
\[
\delta_2(K,K')=\min_{R\in SO(n)}\|h_{RK}-h_{K'}\|_2.
\]
It is well defined since $R\mapsto\|h_{RK}-h_{K'}\|_2$ is easily checked to be continuous with respect to the standard topology of $SO(n)$ (see Section \ref{sec:square} for more details) and $SO(n)$ is compact. In particular, $\varnothing\ne\argmin\{R\mapsto\|h_{RK}-h_{K'}\|_2\}$. For the sequence of symmetrals $(K_m)$, $\delta_2$ satisfies the following property.

\begin{lemma}\label{lem:energy}
Given the sequences $(K_m)\subset\K^n$ and $(f_m)\subset C(\Sn)$ as above, for every $j>i\geq1$
\begin{equation}\label{eq:energy}
\delta_2(K_i,K_j)^2\leq\|f_i\|_2^2-\|f_j\|_2^2.
\end{equation}
\end{lemma}

\begin{proof}
Fix $j>i\geq1$. Repeated use of~\eqref{eq:order} and preservation of pointwise order gives
\[
0\leq f_j\leq T_jf_{j-1}\leq T_j\cdots T_{i+1}f_i=g.
\]
Expanding the operators, we can write
\[
g=\sum_{k=1}^{2^{j-i}}\frac{1}{2^{j-i}}\mathcal{U}_{Q_k}f_i,
\]
where the $Q_k\in O(n)$ are compositions of reflections. Hence
\begin{align*}
\sum_{k=1}^{2^{j-i}}\frac{1}{2^{j-i}}\|\mathcal{U}_{Q_k}f_i-f_j\|_2^2
&=\|f_i\|_2^2+\|f_j\|_2^2-2\langle g,f_j\rangle_2\\
&=\|f_i\|_2^2-\|f_j\|_2^2-2\langle g-f_j,f_j\rangle_2\\
&\leq\|f_i\|_2^2-\|f_j\|_2^2.
\end{align*}
Since this is an average of nonnegative numbers, at least one term satisfies the same upper bound. Thus, for some $k_0$,
\[
\|\mathcal{U}_{Q_{k_0}}f_i-f_j\|_2^2\leq\|f_i\|_2^2-\|f_j\|_2^2.
\]
If $\det Q_{k_0}=-1$, replace $Q_{k_0}$ by $Q_{k_0}\rho_{H_i}\in SO(n)$. This does not change $Q_{k_0}K_i$ because $K_i$ is $H_i$-symmetric. The last display is therefore an admissible competitor for $\delta_2(K_i,K_j)$ and proves~\eqref{eq:energy}.
\end{proof}

Next, observe that for every $A,A',K,K'\in\K^n$
\begin{equation}\label{eq:continuity}
|\delta_2(A,A')-\delta_2(K,K')|\leq\|h_A-h_K\|_2+\|h_{A'}-h_{K'}\|_2.
\end{equation}
Indeed, notice that for $R\in SO(n)$
\[
\|h_{RA}-h_{A'}\|_2\leq\|h_{RA}-h_{RK}\|_2+\|h_{RK}-h_{K'}\|_2+\|h_{A'}-h_{K'}\|_2.
\]
In particular,
\[
\delta_2(A,A')\leq\delta_2(K,K')+\|h_A-h_K\|_2+\|h_{A'}-h_{K'}\|_2,
\]
and~\eqref{eq:continuity} is proved inverting the role of $A,A'$ and $K,K'$.

To finish the proof of Theorem~\ref{thm:one}, by the Blaschke selection theorem~\cite[Theorem~1.8.7]{schneider} there exists a subsequence $(K_{m_j})$ converging to some $L\in\K^n$ in the Hausdorff metric. For every $m \in \N$ apply Lemma~\ref{lem:energy} to the indices $m_j>m$ and let $j\to\infty$ to obtain
\[
\delta_2(K_m,L)^2\leq\|f_m\|_2^2-s^2.
\]
The right-hand side tends to zero. Now, for every $m$, choose a rotation $R_m$ minimizing $\delta_2(K_{m},L)$. Then $h_{R_mK_m}\to h_L$ in $L^2(\Sn)$. Since all these convex bodies lie in the fixed ball $aB^n$, convergence in the $L^2$ support-function metric is equivalent to Hausdorff convergence, see~\cite[Notes for Section~1.8, Note~9]{schneider}. Hence $R_mK_m\to L$.

\subsection{A 2-dimensional example: Behaviour on ellipses}

Before completing the proof of Theorem \ref{thm:main} in the next section, we present an easier proof restricted to centred ellipses in $\R^2$ to give a geometric intuition of \eqref{eq:dichotomy}, which we are going to prove in full generality in the next section. We recall that for every $1$-dimensional subspace $H$, the set $S_HE$ is still an ellipse, and its axes are aligned along $H$ and $H^\perp$.

\begin{theorem}\label{thm:ellipse}
Consider a sequence $(H_m)\subset\Gr_{2,1}$. Then, for every centred ellipse $E_0\subset\R^2$ with semiaxes of lengths $a_0>b_0>0$, the sequence of symmetrals
\[
E_m=S_{H_m}E_{m-1}
\]
has necessarily one of the following mutually exclusive behaviours:
\begin{itemize}
\item The sequence $E_m$ converges to some ellipse $E$.
\item The sequence $E_m$ does not converge, but converges in shape to some ellipse $E$ and, denoting by $\alpha_m\in[0,\pi/2]$ the smaller angle between the major axes of $E_m$ and $E_{m-1}$,
\[
\sum_{m\in\N}\alpha_m=\infty,\qquad
\sum_{m\in\N}\alpha_m^2<\infty.
\]
\end{itemize}
\end{theorem}

\begin{proof}
If $E_m$ converges, then we are in the first case. Suppose otherwise. Theorem~\ref{thm:one} implies convergence in shape. Write $\lambda_m=a_m^2$ and $\mu_m=b_m^2$, where $a_m\geq b_m$ are the semiaxis lengths. We represent a centred ellipse by a positive definite matrix $A$ through $\{x\in\R^2:x^T A^{-1}x\leq 1\}$. In coordinates in which $H$ is the first axis, an ellipse represented by $A=\bigl(\begin{smallmatrix}d&b\\b&c\end{smallmatrix}\bigr)$ is sent by $S_H$ to the ellipse represented by $\diag(d,\det A/d)$. Consequently, $\lambda_m$ is nonincreasing, $\mu_m$ is nondecreasing, and $\lambda_m\mu_m=\lambda_0\mu_0$. Hence $\lambda_m\to\lambda_\infty$ and $\mu_m\to\mu_\infty$, with $\lambda_\infty\geq\mu_\infty$. The inequality is strict, since otherwise $\sqrt{\mu_m}B^2\subseteq E_m\subseteq\sqrt{\lambda_m}B^2$ would force convergence to a disk. No symmetral $E_m$ can itself be a disk, since all subsequent symmetrals would then be the same. Thus, all the major axes are well defined. Finally,
\[
\sum_{m\in\N}\alpha_m=\infty,
\]
because a finite sum would make the sequence of major-axis lines a Cauchy sequence in the angular metric. Together with convergence of the semiaxes, this would imply convergence of $E_m$.

Now, observe that
\begin{equation}\label{eq:ellipse}
\lambda_{m-1}-\lambda_m\geq(\lambda_{m-1}-\mu_{m-1})\sin^2\alpha_m\geq(\lambda_\infty-\mu_\infty)\sin^2\alpha_m.
\end{equation}
Indeed, for two successive symmetrals $E_{m-1}$ and $E_m$, the major axis shifts in one of the two following ways: Either it aligns with $H_m$, in which case
\[
\lambda_m=\lambda_{m-1}-(\lambda_{m-1}-\mu_{m-1})\sin^2\alpha_m,
\]
or with $H_m^\perp$, and
\[
\lambda_m=\frac{\lambda_{m-1}\mu_{m-1}}{\mu_{m-1}+(\lambda_{m-1}-\mu_{m-1})\sin^2\alpha_m}.
\]
In both cases, we obtain~\eqref{eq:ellipse}. To conclude, the elementary inequality $\sin\alpha\geq 2\alpha/\pi$ for $\alpha\in[0,\pi/2]$, together with~\eqref{eq:ellipse}, gives
\[
\alpha_m^2\leq\frac{\pi^2}{4(\lambda_\infty-\mu_\infty)}(\lambda_{m-1}-\lambda_m).
\]
Thus,
\[
\sum_{m\in\N}\alpha_m^2
\leq\frac{\pi^2}{4(\lambda_\infty-\mu_\infty)}
\sum_{m\in\N}(\lambda_{m-1}-\lambda_m)
=\frac{\pi^2}{4(\lambda_\infty-\mu_\infty)}(\lambda_0-\lambda_\infty)<\infty.
\]
\end{proof}

\section{Choosing the right rotations in Theorem \ref{thm:main}}\label{sec:square}

To obtain the full statement of Theorem~\ref{thm:main}, we first work with a single sequence of symmetrals. The main idea is the following: Working on the orbit of the possible limits of the sequence $(K_m)$, one can iteratively choose the rotation $R_{m+1}$ from Theorem~\ref{thm:one} in a suitably small neighbourhood of $R_m$ inside $SO(n)$. To this end, we work on $SO(n)$ as a Lie group and a Riemannian manifold. We start by recalling the corresponding metric. In this section, for ease of notation we write $G=SO(n)$. Then, on $G$ we consider its standard Riemannian metric, which, at the identity, reads as
\[
\langle X,Y\rangle_{\mathfrak{g}}=\frac12\operatorname{tr}(X^TY),\qquad X,Y\in\mathfrak g=\mathfrak{so}(n).
\]
Here, $\mathfrak{so}(n)$ is the Lie algebra of $SO(n)$. See, for example,~\cite{brocker} for this and other basic definitions. As it is standard, the metric on the whole $G$ is defined by left translates. Recall that both left and right translates are isometries for this metric. In particular, for every $A,B,Q\in SO(n)$
\begin{equation}\label{eq:biinvariant}
d_G(AQ,BQ)=d_G(QA,QB)=d_G(A,B),\qquad d_G(QAQ^{-1},I)=d_G(A,I).
\end{equation}

We will need the following lemma. Given a real finite-dimensional inner-product space $V$, $O(V)$ denotes its orthogonal group and $\|\cdot\|_V$ its norm.

\begin{lemma}\label{lem:local}
Consider an orthogonal smooth representation of $G$ on $V$. That is, we are given a smooth group homomorphism
\[
\begin{aligned}
\mathcal{U}:G&\longrightarrow O(V)\\
Q&\longmapsto \mathcal{U}_Q.
\end{aligned}
\]
For a fixed $v\in V$ consider its orbit
\[
M=\{\mathcal{U}_Qv:Q\in G\}.
\]
Then $M$ is a smooth embedded submanifold of $V$. Moreover, there exist $\eta>0$, $C_\eta>0$, such that for every $y,z\in M$ with $\|y-z\|_V<\eta$ there exists $U\in G$ satisfying
\[
\mathcal{U}_Uy=z\qquad\text{and}\qquad d_G(U,I)\leq C_\eta\|y-z\|_V.
\]
\end{lemma}

\begin{proof}
First, suppose that $M=\{v\}$. Then the statement is trivial. Therefore, in the remainder of the proof we suppose that $M$ is not a singleton. Since $G$ is compact, its smooth action on $V$ is proper. See, e.g.,~\cite[Corollary~21.6]{lee}. As orbits of proper Lie group actions are embedded submanifolds, so is $M$, which is compact since it is the continuous image of $G$. This proves the first claim.

To prove the second part observe that the orbit map
\[
\begin{aligned}
p:G&\longrightarrow M\\
Q&\longmapsto \mathcal{U}_Qv
\end{aligned}
\]
is a smooth submersion. By the local section theorem~\cite[Theorem~4.26]{lee} there is a smooth local right inverse $q$ of $p$ with $q(v)=\Id$. Then, upon choosing $\eta>0$ sufficiently small, there exists a constant $C_\eta>0$ such that for every $w\in M$ with $\|v-w\|_V<\eta$ holds
\[
d_G(q(w),\Id)\leq C_\eta\|w-v\|_V.
\]
The conclusion now follows by a standard translation argument. For $y,z\in M$ with $\|y-z\|_V<\eta$, choose $Q\in G$ with $\mathcal{U}_Qv=y$ and set $w=\mathcal{U}_{Q^{-1}}z$. Then $U=Qq(w)Q^{-1}$ satisfies
\[
\mathcal{U}_U y=z,\qquad d_G(U,\Id)=d_G(q(w),\Id)\leq C_\eta\|y-z\|_V,
\]
where bi-invariance of $d_G$ and orthogonality of the representation were used.
\end{proof}

\subsection{Square-summable convergence in shape}

We now prove that the rotations in Theorem \ref{thm:one} can be chosen to have successive square-summable distance in $G$. This behaviour already appeared in \cite{bbgv}, and Bianchi and Gronchi \cite{exposition} suggested the author that this phenomenon should extend to all sequences of iterated Steiner symmetrals of convex bodies.
\begin{theorem}\label{thm:square}
Consider $(K_m)$ as in Theorem~\ref{thm:main}. If $(K_m)$ does not converge, then for any subsequential limit $L$ of $(K_m)$ the rotations $(R_m)\subset G$ giving $R_mK_m\to L$ can be chosen so that
\begin{equation}\label{eq:square}
\sum_{m\in\N}d_G(R_m,R_{m-1})^2<\infty.
\end{equation}
\end{theorem}

\begin{proof}
Retain the notation $a,f_m,T_m,s$ from Section~\ref{sec:shape}, and fix a limit $L$ as in the proof of Theorem~\ref{thm:one}. We start by refining the main argument in the proof of Lemma~\ref{lem:energy}. Since $T_m$ is self-adjoint and $T_mf_m=f_m$,
\[
\langle f_{m-1},f_m\rangle_2=\langle T_mf_{m-1},f_m\rangle_2\geq\|f_m\|_2^2.
\]
Consequently,
\[
\|f_m-f_{m-1}\|_2^2\leq\|f_{m-1}\|_2^2-\|f_m\|_2^2,
\]
and therefore
\begin{equation}\label{eq:steps}
\sum_{m\in\N}\|f_m-f_{m-1}\|_2^2\leq\|f_0\|_2^2-s^2<\infty.
\end{equation}
Thus, $\|f_m-f_{m-1}\|_2\to0$. Write $f=h_L+a$ and consider the compact orbit
\[
\mathcal O=\{\mathcal{U}_Qf:Q\in G\}\subset L^2(\Sn).
\]
By the considerations above, $d_{L^2}(f_m,\mathcal O)^2=\delta_2(K_m,L)^2\to0$. Although $\mathcal O$ is homeomorphic to the finite-dimensional homogeneous space $G/\Stab_G(f)$, its linear span may be infinite-dimensional. In order to be able to apply Lemma~\ref{lem:local}, we now get to a reduction argument.

For $j\in\N\cup\{0\}$ let $\mathcal H_j^n$ be the space of spherical harmonics of degree $j$ and set
\[
V_N=\bigoplus_{j=0}^N\mathcal H_j^n,\qquad P_N:L^2(\Sn)\to V_N.
\]
These orthogonal projections commute with rotations and satisfy $P_Nf\to f$ in $L^2$. Define
\[
G_N=\{Q\in G:\mathcal{U}_QP_Nf=P_Nf\}.
\]
Each $G_N$ is a closed subgroup of $G$. Since the projections commute with rotations and $P_NP_{N+1}=P_N$, these groups form a decreasing sequence. The descending-chain theorem for closed subgroups of a compact Lie group~\cite[p.~136]{brocker} gives $N_0$ such that $G_N=G_{N_0}$ for every $N\geq N_0$. Moreover,
\[
\bigcap_{N=0}^{\infty}G_N=\Stab_G(f),
\]
because for $Q$ in the intersection, $\|\mathcal{U}_Qf-f\|_2\leq2\|f-P_Nf\|_2\to0$. Therefore
\begin{equation}\label{eq:stabilizer}
G_{N_0}=\Stab_G(f).
\end{equation}
Write $P=P_{N_0}$, $V=V_{N_0}$, and $v=Pf$. Equality of stabilizers implies that $P$ is injective on the orbit $\mathcal O$. Let
\[
M=\{\mathcal{U}_Qv:Q\in G\}\subset V,\qquad x_m=Pf_m.
\]
Then $d_V(x_m,M)\to0$ and $\|x_m-x_{m-1}\|_V\leq\|f_m-f_{m-1}\|_2$.

If $M$ is a singleton, injectivity of $P|_{\mathcal O}$ makes $\mathcal O$ a singleton, and $d_{L^2}(f_m,\mathcal O)\to0$ would imply convergence of $(K_m)$, contrary to the hypothesis. Thus $M$ is not a singleton. By Lemma~\ref{lem:local}, $M$ is a compact embedded submanifold of $V$. By the tubular neighbourhood theorem~\cite[Theorem~6.24]{lee}, nearest-point projection $\pi$ onto $M$ is smooth on a sufficiently small tube, and its differential is bounded on a smaller closed tube by compactness. Since $d_V(x_m,M)\to0$ and $\|x_m-x_{m-1}\|_V\to0$, the segment from $x_{m-1}$ to $x_m$ lies in the smaller tube for all sufficiently large $m$. In particular, the points $y_m=\pi(x_m)$ satisfy
\[
\|y_m-y_{m-1}\|_V\leq C_\pi\|x_m-x_{m-1}\|_V.
\]

Choose $m_0$ sufficiently large that Lemma~\ref{lem:local} gives $U_m\in G$ for $m>m_0$ such that
\[
\mathcal{U}_{U_m}y_{m-1}=y_m,\qquad d_G(U_m,\Id)\leq C_\eta\|y_m-y_{m-1}\|_V.
\]
Additionally, consider $Q_{m_0}\in G$ with $\mathcal{U}_{Q_{m_0}}v=y_{m_0}$. Set $R_m=Q_{m_0}^{-1}$ for $0\leq m\leq m_0$, and for $m>m_0$ define
\[
Q_m=U_mQ_{m-1},\qquad R_m=Q_m^{-1}.
\]
Then $\mathcal{U}_{R_m}y_m=v$, and
\begin{align*}
d_G(R_m,R_{m-1})^2
&=d_G(U_m,\Id)^2\leq C_\eta^2\|y_m-y_{m-1}\|_V^2\\
&\leq C_\eta^2 C_\pi^2\|x_m-x_{m-1}\|_V^2\\
&\leq C_\eta^2 C_\pi^2\|f_m-f_{m-1}\|_2^2\\
&\leq C_\eta^2 C_\pi^2\bigl(\|f_{m-1}\|_2^2-\|f_m\|_2^2\bigr).
\end{align*}
Thus,~\eqref{eq:square} follows by~\eqref{eq:steps}.

Finally, we show that $R_mK_m\to L$. Put $g_m=\mathcal{U}_{R_m}f_m$. Since $P$ commutes with rotations,
\begin{equation}\label{eq:projected}
Pg_m-v=\mathcal{U}_{R_m}(x_m-y_m)\longrightarrow0.
\end{equation}
Also,
\[
d_{L^2}(g_m,\mathcal O)=d_{L^2}(f_m,\mathcal O)\longrightarrow0.
\]
Since $\mathcal O$ is compact, every subsequence of $(g_m)$ has a further subsequence converging in $L^2$ to some $g\in\mathcal O$, and~\eqref{eq:projected} implies $Pg=Pf$. Since $P|_{\mathcal O}$ is injective, $g=f$. Hence $g_m\to f$ in $L^2$. Since
\[
g_m=h_{R_mK_m}+a,\qquad f=h_L+a,
\]
the equivalence of the $L^2$ support-function and Hausdorff topologies on uniformly bounded convex bodies~\cite[Notes for Section~1.8, Note~9]{schneider} yields $R_mK_m\to L$.
\end{proof}

\begin{remark}
We use a dimension reduction from $\mathcal O$ to $M$, since even if $\mathcal O$ is finite dimensional, its linear span need not be. E.g.\ support functions of noncircular ellipses have infinitely many Fourier coefficients, so their rotation orbits have infinite-dimensional linear span.
\end{remark}

\subsection{Finding a common sequence of rotations}\label{subsec:common}

We are left to show that the rotations in Theorem~\ref{thm:main} can be chosen independently of the initial body. We start by considering finitely many initial bodies and applying the preceding arguments to their symmetrals simultaneously. For a tuple $\mathbf K=(K_1,\ldots,K_r)$ and $Q\in G$, we write $Q\mathbf K=(QK_1,\ldots,QK_r)$. Convergence of tuples is understood in the Hausdorff metric in each component.

\begin{lemma}\label{lem:finite-family}
Consider a sequence $(H_m)\subset\Gr_{n,n-1}$ and a finite family $K_0^{(1)},\ldots,K_0^{(r)}\in\K^n$, with $r\geq1$. For each $\ell=1,\ldots,r$ let $K_m^{(\ell)}=S_{H_m}K_{m-1}^{(\ell)}$, and write $\mathbf K_m=(K_m^{(1)},\ldots,K_m^{(r)})$. Then any two subsequential limits of $\mathbf K_m$ differ by a common rotation. Moreover, for each subsequential limit $\mathbf L$ there is a sequence $(R_m)\subset G$ satisfying~\eqref{eq:square} such that $R_m\mathbf K_m\to\mathbf L$.
\end{lemma}

\begin{proof}
Choose $a>0$ such that $K_0^{(\ell)}\subseteq aB^n$ for every $\ell$. In $\mathcal H=(L^2(\Sn))^r$, endowed with the product scalar product, put
\[
\mathbf f_m=(h_{K_m^{(1)}}+a,\ldots,h_{K_m^{(r)}}+a),
\qquad E_m=\|\mathbf f_m\|_{\mathcal H}^2.
\]
We let $\mathcal{U}_Q$ and $T_m$ act in each component as in Section~\ref{sec:shape}. By~\eqref{eq:order}, $0\leq\mathbf f_m\leq T_m\mathbf f_{m-1}$ componentwise, so $E_m$ decreases to a limit $E_\infty$. Notice that in the expansion of $T_j\cdots T_{i+1}$ the same compositions of reflections appear in each component. Moreover, if $\det Q=-1$, then $Q\rho_{H_i}\in G$ and $\mathcal{U}_{Q\rho_{H_i}}\mathbf f_i=\mathcal{U}_Q\mathbf f_i$, since each $K_i^{(\ell)}$ is $H_i$-symmetric. Thus the proof of Lemma~\ref{lem:energy} gives, for $j>i\geq1$,
\[
\min_{Q\in G}\|\mathcal{U}_Q\mathbf f_i-\mathbf f_j\|_{\mathcal H}^2\leq E_i-E_j.
\]

Now let $\mathbf L=(L_1,\ldots,L_r)$ be a subsequential limit and put $\mathbf f=(h_{L_1}+a,\ldots,h_{L_r}+a)$. Passing to the corresponding subsequence in the preceding inequality gives
\[
d_{\mathcal H}(\mathbf f_i,\mathcal O)^2\leq E_i-E_\infty,
\qquad \mathcal O=\{\mathcal{U}_Q\mathbf f:Q\in G\}.
\]
In particular, since $\mathcal O$ is compact, every subsequential limit of $\mathbf f_m$ lies in $\mathcal O$. The first assertion follows since support functions determine their bodies.

If $\mathbf K_m$ converges, we can take $R_m=\Id$. Otherwise, notice that $T_m\mathbf f_m=\mathbf f_m$. As in the proof of Theorem~\ref{thm:square}, we then obtain
\[
\|\mathbf f_m-\mathbf f_{m-1}\|_{\mathcal H}^2\leq E_{m-1}-E_m,
\qquad
\sum_{m\in\N}\|\mathbf f_m-\mathbf f_{m-1}\|_{\mathcal H}^2\leq E_0-E_\infty.
\]
Apply $P_N$ in every component. The stabilizers of the resulting projections of $\mathbf f$ form a decreasing sequence with intersection $\Stab_G(\mathbf f)$. By the descending-chain theorem, the projection onto $V_N^r$ is therefore injective on $\mathcal O$ for some $N$. We can now repeat the nearest-point projection and lifting argument of Theorem~\ref{thm:square} in $V_N^r$. It gives rotations satisfying~\eqref{eq:square} and $\mathcal{U}_{R_m}\mathbf f_m\to\mathbf f$ in $\mathcal H$. By the equivalence of the support-function and Hausdorff topologies on uniformly bounded convex bodies, $R_m\mathbf K_m\to\mathbf L$.
\end{proof}

We now choose a finite family for which adjoining another sequence of symmetrals cannot reduce the stabilizer of the limit tuple. For $\mathbf L=(L_1,\ldots,L_r)$ its stabilizer is $\Stab_G(\mathbf L)=\bigcap_{\ell=1}^r\Stab_G(L_\ell)$.

\begin{proposition}\label{prop:common-rotations}
For every sequence $(H_m)\subset\Gr_{n,n-1}$ there exists a sequence $(R_m)\subset G$ satisfying~\eqref{eq:square} such that, for every $K_0\in\K^n$, the symmetrals $K_m=S_{H_m}K_{m-1}$ satisfy $R_mK_m\to L$ in the Hausdorff metric for some $L\in\K^n$.
\end{proposition}

\begin{proof}
Fix $(H_m)$ and consider the collection $\mathcal D$ of the stabilizers of all subsequential limits of tuples as in Lemma~\ref{lem:finite-family}. Since these are closed subgroups of $G$, the descending-chain theorem gives a minimal member $D$ under inclusion. Choose a tuple of symmetrals $\mathbf K_m$ and a subsequential limit $\mathbf L$ with $\Stab_G(\mathbf L)=D$. By Lemma~\ref{lem:finite-family}, we can fix rotations $(R_m)$ satisfying~\eqref{eq:square} such that $R_m\mathbf K_m\to\mathbf L$.

Consider now any other initial body $K_0\in\K^n$ and let $(K_m)$ be its sequence of symmetrals with respect to the fixed family $(H_m)$. By the Blaschke selection theorem, along a subsequence on which $\mathbf K_m\to\mathbf L$ we may assume that $K_m\to L$ for some $L\in\K^n$. We observe that
\[
\Stab_G(\mathbf L,L)=D\cap\Stab_G(L)\in\mathcal D.
\]
By minimality of $D$, this intersection equals $D$. In particular, every $Q\in D$ fixes $L$.

Let $L'$ be a limit of a subsequence of $R_m K_m$. By compactness we may pass to a further subsequence on which
\[
R_m\to R,\qquad (\mathbf K_m,K_m)\to(\mathbf A,B)
\]
for some tuple $\mathbf A$ of convex bodies and a further $B \in \K^n$. Applying Lemma~\ref{lem:finite-family} to the enlarged tuple, we find $Q\in G$ such that $(\mathbf A,B)=Q(\mathbf L,L)$. Since $R\mathbf A=\mathbf L$, we have $RQ\in D$, and consequently $L'=RB=RQL=L$.
The corrected symmetrals are uniformly bounded, so uniqueness of their subsequential limit gives $R_mK_m\to L$. This proves the result since the rotations were fixed before choosing $K_0$.
\end{proof}

\begin{proof}[Completing the proof of Theorem~\ref{thm:main}]
By Proposition~\ref{prop:common-rotations}, there is a common sequence $(R_m)$ satisfying~\eqref{eq:square}. If the symmetrals converge for every initial body, we may take $R_m=\Id$. Suppose instead that there is an initial body for which $(K_m)$ fails to converge. If $\sum_m d_G(R_m,R_{m-1})<\infty$, the rotations are Cauchy and converge to some $R\in G$. But then $R_mK_m\to L$ implies $K_m\to R^{-1}L$, a contradiction. This completes the proof.
\end{proof}

\section{Final remarks and extensions}\label{sec:extensions}

Apart from $H$-symmetry of the symmetral, a careful inspection of the arguments shows that the key properties of Steiner symmetrization which are used to prove the results above are the following:
\begin{enumerate}
\renewcommand{\labelenumi}{\alph{enumi})}
\item For every $K\in\K^n$ such that $\rho_HK=K$, $S_HK=K$,
\item If $K,L\in\K^n$ satisfy $K\subseteq L$ then $S_HK\subseteq S_HL$,
\item For every $K\in\K^n$ the inclusion $S_HK\subseteq M_HK$ holds.
\end{enumerate}
In fact, (a) and (b) were used to keep the symmetrals in a fixed ball, and this also follows directly from (c) since $M_H$ is monotone and fixes $aB^n$. It follows that if we were to identify maps of $\K^n$ into itself producing $H$-symmetric sets and satisfying (c), the proofs in the previous sections would apply to these maps as well. For example, Minkowski symmetrization satisfies the required assumptions. As a matter of fact, Bianchi, Gardner, and Gronchi~\cite[Corollary~7.3]{symgeom} proved that (a) and (b) give (c) under one additional condition. These maps are particular examples of $H$-symmetrizations, as defined in~\cite{symgeom,symproc}. We report a version of~\cite[Corollary~7.3]{symgeom} adapted to our context. Another application of this result to convergence of symmetrization processes can be found in~\cite{ulivelli2023}.

For a linear subspace $H$ of $\R^n$, we use the same notation $\rho_H$ for the orthogonal map that fixes $H$ and acts by $-{\rm Id}$ on $H^\perp$. The Minkowski symmetral $M_HC$ is defined by the same formula as in Section~\ref{sec:shape}. We denote by $F_HK$ the fibre symmetral of $K$, obtained by replacing each non-empty section $K\cap(x+H^\perp)$, $x\in H$, by its Minkowski symmetral within that fibre with respect to reflection about $x$. In particular, $F_HK=S_HK$ when $H$ is a hyperplane.

\begin{proposition}\label{prop:sandwich}
Consider a map
\[
\diamondsuit:\left(\bigcup_{i=0}^{n-1}\Gr_{n,i}\right)\times\K^n\longrightarrow\K^n
\]
such that for every $H$ and $K$ one has $\rho_H\diamondsuit_HK=\diamondsuit_HK$. Suppose, moreover, that the analogues of (a) and (b) hold with $S_H$ replaced by $\diamondsuit_H$, and that
\[
\diamondsuit_H(K+x)=K
\]
whenever $K$ is $H$-symmetric and $x\in H^\perp$. Then
\[
F_HK\subseteq\diamondsuit_HK\subseteq M_HK
\]
for every such $H$ and every $K\in\K^n$.
\end{proposition}

The corresponding convergence result is as follows. As in Theorem~\ref{thm:main}, we allow the rotations to depend on the sequence of subspaces alone.

\begin{corollary}\label{cor:general}
Consider a sequence $(H_m)$ of proper linear subspaces of $\R^n$, possibly of different dimensions. Then there exists a sequence $(R_m)\subset SO(n)$ satisfying~\eqref{eq:square} with the following property. For every $K\in\K^n$ and every sequence of $H_m$-symmetrizations $\diamondsuit_m$ such that $\diamondsuit_m A\subseteq M_{H_m}A$ for every $A\in\K^n$, the sequence
\[
K_0=K,\qquad K_m=\diamondsuit_m K_{m-1},\quad m\in\N,
\]
satisfies
\[
\dH(R_mK_m,L)\to0
\]
for some $L\in\K^n$. In particular, the same rotations apply to fibre and Minkowski symmetrizations and to the symmetrizations in Proposition~\ref{prop:sandwich}. If one of these uncorrected sequences does not converge, the rotations satisfy~\eqref{eq:dichotomy}.
\end{corollary}

\begin{proof}
The containment in the Minkowski symmetral implies~\eqref{eq:order}, and $H_m$-symmetry gives $T_mf_m=f_m$. These observations still hold for a tuple of processes that use different symmetrizations along the same subspaces. When each $H_m$ is a hyperplane, the proof of Theorem~\ref{thm:main} therefore applies without further changes.

For arbitrary subspaces we need to reconsider the last step in the proof of Lemma~\ref{lem:energy}, since $\rho_{H_i}$ can have determinant $1$. Before this step, the proof gives~\eqref{eq:energy} with the minimum over $O(n)$, also for tuples. As in the proof of Theorem~\ref{thm:one}, it follows that $\mathbf f_m$ approaches the $O(n)$-orbit of the support-function tuple of any subsequential limit $\mathbf L$. This orbit consists of at most two $SO(n)$-orbits. If there are two, they are compact and disjoint, hence have positive distance. On the other hand,~\eqref{eq:steps} holds in the product norm and shows that $\|\mathbf f_m-\mathbf f_{m-1}\|\to0$. Thus, for sufficiently large $m$, the tuple cannot pass from a small neighbourhood of one orbit to a small neighbourhood of the other. Since $\mathbf L$ is a subsequential limit, $\mathbf f_m$ approaches the $SO(n)$-orbit corresponding to $\mathbf L$.

We can now use the proof of Theorem~\ref{thm:square} for tuples, followed by the minimal-stabilizer argument in Section~\ref{subsec:common}. In choosing a minimal stabilizer there, we consider the tuple limits of all the processes in the statement, so that the resulting rotations are common to them all. The same argument at the end of Section~\ref{sec:square} gives the first half of~\eqref{eq:dichotomy} if any uncorrected sequence fails to converge.
\end{proof}

We conclude with a further extension that applies to Minkowski symmetrization. As pointed out in~\cite{klain}, the behaviour of Minkowski symmetrization on $\K^n$ can be used to describe the one on the family of non-empty compact subsets of $\R^n$, which we denote by $\C^n$. The corresponding result in this case looks as follows.

\begin{corollary}\label{cor:compact}
For $C\in\C^n$ and a sequence $(H_m)$ as in Corollary~\ref{cor:general}, consider the sequence
\[
C_0=C,\qquad C_m=M_{H_m}C_{m-1},\quad m\in\N.
\]
Then the same rotations $(R_m)\subset SO(n)$ as in Corollary~\ref{cor:general} give
\[
\dH(R_mC_m,L)\to0
\]
for some $L\in\K^n$. In particular, one sequence of rotations works for every $C\in\C^n$. If $(C_m)$ does not converge for some $C\in\C^n$, these rotations satisfy~\eqref{eq:dichotomy}.
\end{corollary}

\begin{proof}
Put $K_m=\conv C_m$. First, recall that the convex hull commutes with Minkowski addition and, therefore, with Minkowski symmetrization. In particular, $(K_m)$ is the Minkowski symmetrization process starting from $\conv C$. Moreover, the estimate from the Shapley--Folkman--Starr theorem \cite[Theorem 3.1.6]{schneider} gives
\[
\dH(C_m,K_m)\leq\frac{\sqrt n}{2^m}\diam C\longrightarrow0.
\]
See~\cite[Theorem~3.4 and Corollary~3.5]{klain} for a similar argument. Applying the Corollary \ref{cor:general} to $(K_m)$ and using rotation invariance of $\dH$, we obtain convergence of $(R_mC_m)$ to the same limit as $(R_mK_m)$. Finally, the estimate shows that $(C_m)$ converges if and only if $(K_m)$ converges, so the last claim also follows Corollary \ref{cor:general}.
\end{proof}

\paragraph{Acknowledgments.} The author was supported, in part, by the Austrian Science Fund (FWF): 10.55776/PAT3787224. The author thanks Georg Hofst\"atter for insightful discussions.

\medskip
\paragraph{AI disclosure.} This manuscript, both in writing and mathematical proofs, was prepared with the assistance of AI. In particular, GPT-6 Astra was used. The mathematical results were assisted and refined through various rounds of prompting. The final form of proofs and writing in the manuscript were handcrafted by the author, and double-checked with AI for typos and stylistic suggestions.

\medskip
\noindent Jacopo Ulivelli\\
Institut f\"ur Diskrete Mathematik und Geometrie\\
TU Wien\\
Wiedner Hauptstra\ss e 8--10/1046\\
1040 Wien, Austria\\
\textit{e-mail}: \texttt{jacopo.ulivelli@tuwien.ac.at}

\end{document}